\documentclass[a4paper,11pt]{amsart}

\newtheorem{theorem}{Theorem}[section]
\newtheorem{lemma}[theorem]{Lemma}
\newtheorem{proposition}[theorem]{Proposition}

\numberwithin{equation}{section}

\newcommand{\wh}{\widehat}
\newcommand{\Up}{\Upsilon}
\newcommand{\SL}{{\operatorname{SL}}}
\newcommand{\SO}{{\operatorname{SO}}}
\newcommand{\nd}{\nabla}

\newcommand{\ce}{{\mathcal E}}
\newcommand{\cq}{{\mathcal Q}}

\newcommand{\bg}{\mbox{\boldmath{$g$}}}
\newcommand{\bR}{\mathbb{R}}
\newcommand{\bC}{\mathbb{C}}
\newcommand{\bN}{\mathbb{N}}

\newcommand{\Rho}{P}
\newcommand{\RhoP}{{\mbox{\sf P}}}
\newcommand{\RP}{{\mbox{\sf R}}}
\newcommand{\CP}{{\mbox{\sf C}}}

\def\sideremark#1{%
\ifvmode\leavevmode\fi\vadjust{\vbox to0pt{\vss%
 \hbox to 0pt{\hskip\hsize\hskip1em%
 \vbox{\hsize3cm\tiny\raggedright\pretolerance10000%
 \noindent #1\hfill}\hss}\vbox to8pt{\vfil}\vss}}}

\begin{document}

\title{
Conformally invariant powers of the Laplacian --- 
A complete non-existence theorem}

\begin{abstract}
  We show that on conformal manifolds of even dimension $n\geq 4$
  there is no
  conformally invariant natural differential operator between density
  bundles with leading part a power of the Laplacian $\Delta^{k}$ for
  $k>n/2$.  This shows that a large class of invariant operators on
  conformally flat manifolds do not generalise to arbitrarily curved
  manifolds and that the theorem of Graham, Jenne, Mason and Sparling,
  asserting the existence of curved version of $\Delta^k$ for $1\le
  k\le n/2$, is sharp.
\end{abstract}

\thanks{ARG gratefully acknowledges support from the Royal Society of 
New Zealand via a Marsden Grant (grant no.\ 02-UOA-108).
KH gratefully acknowledges support from the Japan Society for the Promotion of Science.}



\author{A. Rod Gover}
\date{}
\address{Department of Mathematics\\
  The University of Auckland\\
  Private Bag 92019\\
  Auckland 1\\
  New Zealand}\email{gover@math.auckland.ac.nz}
  
\author{Kengo Hirachi}
\address{
Graduate School of Mathematical Sciences\\ 
University of Tokyo\\ 
3-8-1 Komaba, Megro, Tokyo 153-8914\\ 
Japan}\email{hirachi@ms.u-tokyo.ac.jp}

\maketitle

\pagestyle{myheadings}
\markboth{A.R. GOVER AND K. HIRACHI}{Non-existence of conformal Laplacians}

\section{Introduction}

Conformally invariant operators and the equations they determine play
a central role in the study of manifolds with pseudo-Riemannian,
Riemannian, conformal and related structures. This observation dates
back to at least the very early part of the last century when it was shown
that the equations of massless particles on curved space-time exhibit
conformal invariance. In this setting a key operator is the
conformally invariant wave operator which has leading term a
pseudo-Laplacian. The Riemannian signature variant of this operator is
a fundamental tool in the Yamabe problem on compact manifolds. Here
one seeks to find a metric, from a given conformal class, that has
constant scalar curvature. Recently it has become clear that higher
order analogues of these operators, viz.\ conformally invariant
operators on weighted functions (i.e., conformal densities) with
leading term a power of the Laplacian, have a central role in
generating and solving other curvature prescription problems as well
as other problems in geometric spectral theory and mathematical
physics \cite{Brsharp,CY,GZ}.

In the flat setting, the existence of such operators dates back to
\cite{JV}, where it is shown that, on 4-dimensional Minkowski space,
for $k\in\bN=\{1,2,\dots\}$, the $k^\text{th}$ power of the
flat wave operator $\Delta^k$,
acting on densities
of the appropriate weight, is invariant under the action of the
conformal group.  More generally, if $\ce[w]$ denotes the space of
conformal densities of weight $w\in\bR$, then on a flat conformal
manifold of dimension $n\geq 3$ (and
any signature) there exists, for each $k\in\bN$, a unique 
conformally invariant operator
$$
 \Box_{2k}\colon
 \ce[k-n/2]\to\ce[-k-n/2]
 $$
and the leading part of $\Box_{2k}$ is $\Delta^k$.
Furthermore, this set of operators is complete 
in the sense that it contains all natural conformally invariant
 differential operators (see section 2) between densities.
These facts are easily recovered from the general results in \cite{ESlo}
and references therein.

Many of these operators can be generalised to curved conformal 
manifolds; Graham, Jenne,  Mason and Sparling \cite{GJMS} constructed 
natural conformally invariant operators
$$
P_{2k}\colon\ce[k-n/2]\to\ce[-k-n/2]
$$
with leading term $\Delta^k$ for all $k\in\bN$ and $n\ge3$ 
except for the cases of $n$ even and $k>n/2$. (See the references in 
\cite{GJMS} for earlier constructions of some low order examples.)
They also conjectured that their result is sharp, based partially on 
the fact, proved by Graham \cite{Grno}, that $\Box_6$ in dimension 4
does not admit curved analogue.  Recently \cite{W} added weight to 
this conjecture by establishing the non-existence of a curved analogue
for $\Box_8$ in dimension 6.  In this paper we prove the conjecture. 
We state this as a theorem.

\medskip

\noindent{\bf Theorem. }{\em
If $n\ge4$ is even and $k>n/2$, there is no conformally invariant
natural differential operator between densities with the same
principal part as $\Delta^{k}$.  } 
 
\medskip 

In \cite{Grno} Graham explains that ``the basic reason for the 
non-existence of an invariant curved modification of $\Delta^3$ in 
dimension 4 is the conformal invariance of the classical Bach tensor.''
An analogue of this reasoning still holds true for the proof of our
theorem, although the proof is completely different from that of Graham.
In higher even dimensions we replace the Bach tensor by its analogue,
the Fefferman-Graham obstruction tensor $B_{ab}$, which arises in the
ambient metric construction of \cite{FG2}; see \eqref{bach} in the next
section.  Our strategy for the proof is to construct a curvature
expression that is shown to be non-zero for a class of conformal
metrics for which $B_{ab}=0$, while it is also shown to vanish for the
same class of metrics under the assumption of the existence of the
curved analogue of $\Box_{2k}$. This is a contradiction.  The former
is done by a direct computation using the tractor calculus, which will
be review in section 2; the latter is a consequence of some classical
invariant theory. 

This explanation is somewhat of a simplification. Nevertheless the
proof of the theorem in section \ref{proof} can be viewed as a careful
elaboration of this idea. The proof is greatly simplified by the
use of a special class of metrics and section \ref{metric} is
concerned with showing that this class is non-trivial.

Finally we should point out that there are many other settings where
similar non-existence issues remain to be resolved. In \cite{ESlo}
Eastwood and Slov\'{a}k use and develop some semiholonomic Verma
module theory to prove that in odd dimensions every conformally
invariant operator between irreducible bundles on (locally)
conformally flat manifolds (including spin manifolds) has a curved
analogue. In even dimensions they show that the same is true, save for
an exceptional class of operators.  The class consists of the
operators corresponding dually to those nonstandard nonsingular
homomorphisms which go between the generalised Verma modules at either
extreme of generalised Bernstein-Gelfand-Gelfand resolutions.  This
includes the $\Box_{2k}$ for $k\geq n/2$, as discussed above, but also
many other operators. Some operators in the exceptional class do have
curved analogues. In particular the $\Box_{2k}$ for $k=n/2$.  However
we suspect that otherwise the result of Eastwood and Slov\'{a}k is
sharp.  Similar questions can be asked for many other similar
geometries such as CR structures. In \cite{GoGr} there is a
construction of CR invariant powers of the sub-Laplacian that
generates curved analogues for most, but not all the invariant
operators from the CR flat setting. Once again there is the question
of whether this result is sharp. For more general CR operators the
existence theory is much less developed than in the conformal case.

It is a pleasure to thank Tom Branson and Josef \v{S}ilhan for helpful
discussions. The authors are also appreciative of the referee's perceptive comments.

\section{Conformal geometry and tractor calculus}\label{tractorsect}

We collect here the minimal background materials from conformal
geometry and tractor calculus as required for the proof of the 
theorem. 
The initial development of tractor calculus in conformal
geometry dates back to the work of T.Y. Thomas \cite{T} and was
reformulated and further developed in a modern setting in \cite{BEGo}.
It is intimately related to the Cartan conformal connection; for a 
comprehensive treatment exposing this connection and relating the 
conformal case to the wider setting of parabolic structures see 
\cite{Cap-Gover,Cap-Gover2}.  The calculational techniques, conventions
and notation used here follow \cite{GP} and \cite{Goadv}.

Let $(M,[g])$ be a conformal manifold of dimension $n\ge 3$ and of
signature $(p,q)$. A conformal structure is equivalent to a ray
subbundle $\cq$ of $S^2T^*M$; points of $\cq$ are pairs $(g_x,x)$
where $x\in M$ and $g_x$ is a metric at $x$, each section of $\cq$
gives a metric $g$ on $M$ and the metrics from different sections
agree up to multiplication by a positive function.  The bundle $\cq$
is a principal bundle with group $\bR_+$, and we denote by $\ce[w]$
the vector bundle induced from the representation of $\bR_+$ on $\bR$
given by $t\mapsto t^{-w/2}$.  Sections of $\ce[w]$ are called a {\em
conformal densities of weight $w$} and may be identified with
functions on $\cq$ that are homogeneous of degree $w$, i.e., $f(s^2
g_x,x)=s^w f(g_x,x)$ for any $s\in \bR_+$.  We will often use the same
notation $\ce[w]$ for the space of sections of the bundle. Note that
for each choice of a metric $g$ (i.e., section of $\cq$, which we term
a {\em choice of conformal scale}), we may identify a section $f\in
\ce[w]$ with a function $f_g$ on $M$ by $f_g(x)=f(g_x,x)$.  In
particular, $\ce[0]$ is canonically identified with $C^\infty(M)$.
Finally we emphasise that for $w\ne0$ the bundle $\ce[w]$, by
its definition, depends on the conformal structure.

The operators of our main interest are defined as maps between
densities $P\colon\ce[w]\to\ce[w']$. For each choice of a scale
$g\in[g]$, $P$ induces a map $P_g:C^\infty(M)\to C^\infty(M)$ via the
identifications $\ce[w]\cong C^\infty(M)$. We say that $P$ is a {\em
  natural differential operator} if $P_g$ can be written as a
universal polynomial in covariant derivatives with coefficients
depending polynomially on the metric, its inverse, the curvature
tensor and its covariant derivatives.  The coefficients of natural
operators are called {\em natural tensors}. In the case that they are
scalar they are often also called {\em Riemannian invariants}.  We say $P$
is a {\em conformally invariant differential operator} if it is a
natural operator in this way and is well defined on conformal
structures (i.e.\ is independent of a choice of conformal scale).

We embrace Penrose's abstract index notation \cite{ot} throughout this
paper and indices should be assumed abstract unless otherwise
indicated.  We write $\ce^a$ to denote the tangent bundle on $M$, and
$\ce_a$ the cotangent bundle. We use the notation
$\ce_a[w]=\ce_a\otimes\ce[w]$,
$\ce_{ab}[w]=\ce_a\otimes\ce_b\otimes\ce[w]$ and so on.  An index
which appears twice, once raised and once lowered, indicates a
contraction.  
The symmetric tensor products of cotangent bundle will be written
as
$\ce_{(ab\cdots c)}$ and $\ce_{(ab\cdots c)_0}$ indicates the
completely trace-free subbundle. Similarly, $\ce_{[ab\cdots c]}$
means the skew tensor product, that is, the bundle of
differential forms. We also use these notation to indicate
the projection onto these bundles, e.g.\ $2T_{[ab]}=T_{ab}-T_{ba}$.
These conventions will be extended in an obvious way to the tractor
bundles described below.

Note that there is a tautological function $\bg$ on $\cq$ taking
values in $\ce_{(ab)}$. 
It is the function which assigns to the point
$(g_x,x)\in \cq$ the metric $g_x$ at $x$.  This is homogeneous of
degree 2 since $\bg (s^2 g_x,x) =s^2 g_x$. If $\xi$ is any positive
function on $\cq$ homogeneous of degree $-2$ then $\xi \bg$ is
independent of the action of $\bR_+$ on the fibres of $\cq$, and
so $\xi \bg$ descends to give a metric from the conformal class. Thus
$\bg$ determines and is equivalent to a canonical section of
$\ce_{ab}[2]$ (called the conformal metric) that we also denote $\bg$
(or $\bg_{ab}$). This in turn determines a canonical section $\bg^{ab}$ 
(or $\bg^{-1}$) of $\ce^{ab}[-2]$ with the property that
$\bg_{ab}\bg^{bc}=\delta_a^c$ (where $\delta_a{}^c $ is kronecker
delta, i.e., the section of $\ce^c_a$ corresponding to the identity
endomorphism of the the tangent bundle).  The conformal metric (and
its inverse $\bg^{ab}$) will be used to raise and lower indices.
Given a choice of metric $ g\in[g]$, we write $ \nabla_a$ for the
corresponding Levi-Civita connection. For each choice of metric there
is also a canonical connection on $\ce[w]$ determined by the 
identification of $\ce[w]$ with $C^\infty(M)$, as described above, 
and the exterior derivative on functions. We will also 
call this the Levi-Civita connection and thus for tensors with weight, 
e.g.\ $v_{a}\in\ce_{a}[w]$, there is a connection given by the Leibniz 
rule. With these conventions the
Laplacian $ \Delta$ is given by $\Delta=\bg^{ab}\nd_a\nd_b=
\nd^b\nd_b\,$.  

The Riemannian curvature $R_{ab}{}^c{}_d$, determined by
$$
(\nd_a\nd_b-\nd_b\nd_a)v^c=R_{ab}{}^c{}_dv^d ,
\quad\text{
where}\ v^c\in \ce^c,
$$
can be decomposed into the totally trace-free {\em Weyl curvature}
$C_{abcd}$ and the symmetric {\em
Schouten tensor} $\Rho_{ab}$ according to
\begin{equation} \label{R-PC} 
R_{abcd}=C_{abcd}+2\bg_{c[a}\Rho_{b]d}+2\bg_{d[b}\Rho_{a]c}.
\end{equation}
This defines $P_{ab}$ as a trace modification of the Ricci tensor 
$R_{ab}=R_{ca}{}^c{}_b$:
$$
R_{ab}=(n-2)\Rho_{ab}+\Rho_{c}{}^c\bg_{ab}.
$$
Note that the Weyl tensor has the symmetries 
\begin{equation}\label{C-sym}
 C_{abcd}=C_{[ab][cd]}=C_{cdab},\quad  C_{[abc]d}=0.
\end{equation}
Moreover, it follows from the Bianchi identity that
\begin{equation} \label{Bian} 
 \nd^c C_{abcd} = 2 (n-3)\nd_{[a}\Rho_{b]d}
\end{equation} 
and
\begin{equation} \label{Bian-Weyl} 
 (n-3)\nd_{[a }C_{bc]de} =\bg_{d[a}\nd^sC_{bc]se}-\bg_{e[a}\nd^sC_{bc]sd}.
\end{equation}

Under a {\em conformal transformation}, we replace our
choice of metric $ g$ by the metric $ \wh{g}=e^{2\Up} g$, where
$\Up$ is a smooth function.  The Levi-Civita connection
then transforms as follows:
\begin{equation}\label{transform}
\wh{\nd_a u_b}=\nd_a u_b -\Up_a u_b-\Up_b u_a +\bg_{ab} \Up^c u_c,
 \quad
 \wh{\nd_a \sigma} = \nd_a \sigma +w\Up_a \sigma.
\end{equation}
Here $ u_b\in \ce_b$, $ \sigma\in \ce[w]$, and 
$\Upsilon_a=\nd_a\Up$. The Weyl curvature is
conformally invariant, that is $\widehat{C}=C$, and the
Schouten tensor transforms by
\begin{equation}\label{Rhotrans}
  \widehat{\Rho}_{ab}=\Rho_{ab}-\nd_a \Up_b +O(\Up^2),
\end{equation}
where $O(\Up^2)$ denotes non-linear terms in $\Up$. 

We define $\Rho^{(\ell)}\in\ce_{(a_1\cdots a_{\ell})}$
for $\ell\ge2$
by
$$
 \Rho^{(\ell)}=\Rho_{a_{\ell}\cdots a_1}
 :=\nd_{(a_{\ell}}\cdots \nd_{a_3}\Rho_{a_2 a_1)}.
$$ 
From \eqref{R-PC} and \eqref{Bian} it follows easily that if $n>3$ then the jets
 of $R$ at $p$ can be expressed in terms of $\Rho^{(\ell)}$ and the
 jets of $C$ at $p$.  Note that we can always choose, for each point
 $p\in M$, a representative $g$ from a conformal class such that
\begin{equation}\label{normal}
 \Rho^{(\ell)}(p)=0, \quad \ell\geq 2;
\end{equation}
following \cite{Goadv} we call $g$ a {\em normal scale}.
This is an easy consequence of the conformal variational formula:
\begin{equation}\label{Rhotrans2}
\widehat\Rho_{a_\ell\cdots a_1}
=\Rho_{a_\ell\cdots a_1}
-\nd_{(a_\ell} \cdots\nd_{a_1)}\Up+ O(\Up^2),
\end{equation}
since the terms in $O(\Up^2)$ involve at most $\ell-1$ derivatives of
$\Up$.  In a normal scale, the jets of $R$ at $p$ can be expressed in
terms of the Weyl curvature $C$ and its covariant derivatives at
$p$. 

In dimension $4$, it is well-known that 
$$
 B_{ab}=\nd^d\nd^c C_{cadb}+P^{dc}C_{cadb}
$$ 
is a conformally invariant tensor, called the {\em
Bach tensor}.  The existence of a natural conformally invariant
tensor, taking values in $\ce_{(ab)_0}[2-n]$ and which generalises the
Bach tensor to even dimensions, is deduced in \cite{FG2} where it
arises as the obstruction to the existence of a formal power series
solution to their ambient metric construction. 
We will also denote this {\em Fefferman-Graham
obstruction tensor} by $B_{ab}$. 
While no general explicit expression for
$B_{ab}$ has been given it is easily shown from its origins as an
obstruction that it contains linear terms when we consider
perturbations from the flat metric. Using this, its naturality and
conformal invariance as well as the symmetries and identities
satisfied by the Weyl curvature it is straightforward to deduce that
its linear in curvature term is given (up to non-zero constant
multiple) by
\begin{equation}\label{bach}
 \Delta^{n/2-2}\nd^{c}\nd^dC_{cadb}.
\end{equation}

We next define the standard tractor bundle over $(M,[g])$.
It is a vector bundle of rank $n+2$ defined for each $g\in[g]$
by  $[\ce^A]_g=\ce[1]\oplus\ce_a[1]\oplus\ce[-1]$. 
If $\wh g=e^{2\Up}g$, we identify  
 $(\sigma,\mu_a,\tau)\in[\ce^A]_g$ with
$(\wh\sigma,\wh\mu_a,\wh\tau)\in[\ce^A]_{\wh g}$
by the transformation
\begin{equation}\label{transf-tractor}
 \begin{pmatrix}
 \wh\sigma\\ \wh\mu_a\\ \wh\tau
 \end{pmatrix}=
 \begin{pmatrix}
 1 & 0& 0\\
 \Up_a&\delta_a{}^b&0\\
- \tfrac{1}{2}\Up_c\Up^c &-\Up^b& 1
 \end{pmatrix} 
 \begin{pmatrix}
 \sigma\\ \mu_b\\ \tau
 \end{pmatrix} .
\end{equation}
It is straightforward to verify that these identifications are
consistent upon changing to a third metric from the conformal class,
and so taking the quotient by this equivalence relation defines the
{\em standard tractor bundle} $\ce^A$ over the conformal manifold.
(Alternatively the standard tractor bundle may be constructed as a
canonical quotient of a certain 2-jet bundle or as an associated
bundle to the normal conformal Cartan bundle \cite{Cap-Gover2}.)  The
bundle $\ce^A$ admits an invariant metric $ h_{AB}$ of signature
$(p+1,q+1)$ and an invariant connection, which we shall also denote by
$\nabla_a$, preserving $h_{AB}$.  In a conformal scale $g$, these are
given by
$$
 h_{AB}=\begin{pmatrix}
 0 & 0& 1\\
 0&\bg_{ab}&0\\
1 & 0 & 0
 \end{pmatrix}
\text{ and }
\nabla_a\begin{pmatrix}
 \sigma\\ \mu_b\\ \tau
 \end{pmatrix}
 =
\begin{pmatrix}
 \nabla_a \sigma-\mu_a \\
 \nabla_a \mu_b+ \bg_{ab} \tau +\Rho_{ab}\sigma \\
 \nabla_a \tau - \Rho_{ab}\mu^b  \end{pmatrix}. 
$$
It is readily verified that both of these are conformally well-defined,
i.e., independent of the choice of a metric $g\in [g]$.  Note that
$h_{AB}$ defines a section of $\ce_{AB}=\ce_A\otimes\ce_B$, where
$\ce_A$ is the dual bundle of $\ce^A$. Hence we may use $h_{AB}$ and
its inverse $h^{AB}$ to raise or lower indices of $\ce_A$, $\ce^A$ and
their tensor products.

In computations, it is often useful to introduce 
the `projectors' from $\ce^A$ to
the components $\ce[1]$, $\ce_a[1]$ and $\ce[-1]$ which are determined
by a choice of scale.
They are respectively denoted by $X_A\in\ce_A[1]$, 
$Z_{Aa}\in\ce_{Aa}[1]$ and $Y_A\in\ce_A[-1]$, where
 $\ce_{Aa}[w]=\ce_A\otimes\ce_a\otimes\ce[w]$, etc.
 Using the metrics $h_{AB}$ and $\bg_{ab}$ to raise indices,
we define $X^A, Z^{Aa}, Y^A$. Then we
immediately see that 
$$
Y_AX^A=1,\ \ Z_{Ab}Z^A{}_c=\bg_{bc}
$$
and that all other quadratic combinations that contract the tractor
index vanish. This is summarised in figure~\ref{XYZfigure}. 
\begin{figure}
$$
\begin{array}{l|ccc}
& Y^A & Z^{Ac} & X^{A}
\\
\hline
Y_{A} & 0 & 0 & 1
\\
Z_{Ab} & 0 & \delta_{b}{}^{c} & 0
\\
X_{A} & 1 & 0 & 0
\end{array}
$$
\caption{Tractor inner product}
\label{XYZfigure}
\end{figure}

It is clear from \eqref{transf-tractor} that the first component 
$\sigma$ is independent of the choice of a representative $g$ and 
hence $X^A$ is conformally invariant. 
For $Z^{Aa}$ and $Y^A$, we have the transformation laws:
  \begin{equation}\label{XYZtrans}
  \wh Z^{Aa}=Z^{Aa}+\Up^aX^A, \quad
  \wh Y^A=Y^A-\Up_aZ^{Aa}+O(\Up^2).
\end{equation}

Given a  choice
of conformal  scale we have the corresponding Levi-Civita connection
on tensor and density bundles.  In this setting we can use the coupled 
Levi-Civita tractor connection to act on sections of the tensor product 
of a tensor bundle with a tractor bundle. This is defined by the Leibniz 
rule in the usual way.  For example if
$ u^b V^C \sigma\in \ce^b\otimes \ce^C\otimes \ce[w]=:\ce^{bC}[w]$
then 
$\nd_a u^b V^C \sigma = (\nd_a u^b) V^C \sigma +
u^b(\nd_a V^C) \sigma + u^b V^C \nd_a \sigma$.
Here $\nd$ means the Levi-Civita
connection on $ u^b\in \ce^b$ and $ \sigma\in \ce[w]$,
while it denotes the tractor
connection on $ V^C\in \ce^C$. In particular with this convention we have 
\begin{equation}\label{connids}
\nd_aX_A=Z_{Aa},\quad
\nd_aZ_{Ab}=-\Rho_{ab}X_A-Y_A\bg_{ab}, 
\quad \nd_aY_A=\Rho_{ab}Z_A{}^b.
\end{equation}

Note that if $V$ is a section of $ \ce_{A_1\cdots A_\ell}[w]$,
then the coupled Levi-Civita tractor
connection on $V$ is not conformally invariant but transforms just as the
Levi-Civita connection transforms on densities of the same weight:
$\wh{\nd}_a V = \nd_a V + w\Up_a V$.

Given a choice of conformal scale, the {\em tractor-$D$ operator} 
$$
D_A\colon\ce_{A_1\cdots A_\ell}[w]\to\ce_{AA_1\cdots A_\ell}[w-1]
$$
is defined by 
\begin{equation}\label{Dform}
D_A V:=(n+2w-2)w Y_A V+ (n+2w-2)Z_{Aa}\nabla^a V -X_A\Box V, 
\end{equation} 
 where $\Box V :=\Delta V+w \Rho_b{}^b V$.  This also turns out to be
 conformally invariant as can be checked directly using the formulae
 above (or alternatively there are conformally invariant constructions
 of $D$, see e.g.\ \cite{Gosrni}).

The curvature $ \Omega$ of the tractor connection 
is defined by 
\begin{equation}\label{curvature}
[\nd_a,\nd_b] V^C= \Omega_{ab}{}^C{}_EV^E 
\end{equation}
for $ V^C\in\ce^C$.  Using
\eqref{connids} and the usual formulae for the curvature of the
Levi-Civita connection we calculate (cf. \cite{BEGo})
\begin{equation}\label{tractcurv}
\Omega_{abCE}= Z_C{}^cZ_E{}^e C_{abce}-\frac{2}{n-3}
X_{[C}Z_{E]}{}^e\nd^dC_{abde}. 
\end{equation}
Here, and in the remainder of this section, to simplify the formulae
we have assumed $n\geq 4$.  Since our later discussions are all set in
even dimensions there is no need here for the results in dimension 3.  We
also set 
$$
\Omega_{ABCE}=Z_A{}^aZ_B{}^b\Omega_{abCE},\quad
\Omega_{BsCE}=Z_B{}^b\Omega_{bsCE}.
$$

It is easily verified that $[D_{B},D_{C}]$ vanishes on densities.
For tractors $V\in \ce_{A_1A_2\cdots A_k}[w]$,
it is straightforward to use \eqref{curvature} and \eqref{connids} to show
\begin{equation}\label{Dcomm}
[D_{B},D_{C}] V_{A_1A_2\cdots A_k}=
E_{BCA_1}{}^Q V_{QA_2\cdots F} +
\cdots
+ E_{BCA_k}{}^Q V_{A_1\cdots A_{k-1}Q},
\end{equation} 
where  
\begin{equation}\label{def-E}
\begin{aligned}
E_{ABCE}=(n+2w-2)\Big(&(n+2w-4)\Omega_{ABCE} \\
-2X_{[A}Z_{B}{}_{]}^b&\nd^p \Omega_{pbCE}
+4X_{[A}\Omega_{B]}{}^{s}{}_{CE}\nd_s\Big).
\end{aligned}
\end{equation} 
For our forthcoming calculations, we need to express $Y^AY^CE_{ABCE}$ 
in terms of $C$. The first term of $E$ is killed by contraction with 
$Y^A$ and the last term gives
$$
 4Y^AY^CX_{[A}\Omega_{B]p}{}_{CE}\nd^p
=\frac{-2}{n-3} Z_B{}^b Z_E{}^e(\nd^aC_{aebc})\nd^c.
$$
For the middle term, using \eqref{connids}, we have
$$
-2Y^AY^CX_{[A}Z_{B}{}_{]}^b\nd^p \Omega_{pbCE}
=\frac{1}{n-3}
Z_B{}^bZ_E{}^e\nd^a\nd^cC_{abce} +O(R^2).
$$
Where the $O(R^2)$ indicates 
 non-linear terms in the curvature.
Thus using \eqref{C-sym}, we get
 \begin{equation}\label{YYE}
\begin{aligned}
Y^AY^C&E_{ABCE}=
\frac{n+2w-2}{n-3} \times \\
& Z_{B}{}^b Z_E{}^{e}
\Big(\nd^c\nd^d C_{cedb}+2(\nd^cC_{cedb})\nd^d
\Big) +O(R^2).
\end{aligned}
\end{equation}

\section{Proof of the theorem} \label{proof}

\noindent For the remainder of the paper we restrict to manifolds of even dimension $n$.

Our key tool for the proof is the natural differential 
operator 
$$
 L_g:
 \ce[k-n/2]\to\ce[-k-n/2]
$$
defined by
$$
L_gf:=Y^{A_{k}}\cdots Y^{A_{1}} D_{A_{k}}\cdots D_{A_{1}}f.
$$
 This is
a composition of the conformally invariant operator
$$
D_{A_k}\cdots D_{A_{2}} D_{A_{1}}:\ce[k-n/2]\to
\ce_{A_k\cdots A_2A_{1}}[-n/2] 
$$ 
and the projector, determined by $g$, 
$$
Y^{A_{k}}\cdots Y^{A_{1}}:\ce_{A_k\cdots A_{1}}[-n/2] \to\ce[-k-n/2].
$$
It follows easily from \eqref{Dform} and the relations in figure
\ref{XYZfigure} that $L_g$ has leading part $(-\Delta)^k$.  In view of
(\ref{XYZtrans}) we do not expect $L_g$ to exhibit invariance under
conformal rescaling.  However if $g$ is conformally flat, it turns out 
that $L_g$ is the unique conformally invariant operator between
densities whose leading part is $(-\Delta)^k$; see e.g.  Proposition
2.1 of \cite{GP} or \cite{Gosrni}. 

Our strategy for proving the theorem is as follows.
We study the dependence of $L_gf$ on deformations from
the flat metric $g_0$. 
For a smooth family of
Riemannian (or pseudo-Riemannian) metrics 
$\{g_t\}_{t\in\bR}$ such that $g_0$ is the 
flat metric and $\Up\in C^\infty (M\times\bR)$, let
$$
L[s,t]:=L_{g[s,t]}f(p),\quad\text{where}\quad
g[s,t](x):=e^{2s\Up(x,t)}g_t(x),
$$ 
for each $s\in \bR$.  Here we fix $f\in \ce[k-n/2]$ on
$(M,[g_0])$ and then regard it as a density in $\ce[k-n/2]$ for each
conformal structure $[g_t]$ by the identification $f_{g_t}=f_{g_0}$.
With a view to a contradiction we suppose that there exists a natural
conformally invariant operator $P_{2k}$ between density bundles with leading
term $\Delta^k$, where $k>n/2$. Such an operator in particular gives
an operator on conformally flat spaces and so must appear in
the classification of such operators described in the introduction. 
Thus we have
\begin{equation}\label{nop}
P_{2k}:\ce[k-n/2]\to \ce[-k-n/2].
\end{equation}
(In particular, if  $g$ is conformally flat,
then we have
$(-1)^kP_{2k}=L_g$.)  Then we choose $p\in M$,  and set
$$
 P[t]:=(-1)^k P_{2k}f(p)
$$ 
in the metric $g[s,t]$. Note that, since $P_{2k}$ is conformally
invariant, the right-hand side is independent of $s$.  We compute
$\partial_t\partial_s L[0,0]$ by two methods, which give different
answers.  With some assumptions on the family $g[s,t]$ and on
$f$, we show, by one set of calculations, that 
$$
\partial_t\partial_s L[0,0]\ne0.
$$ 
On the other hand, with the other approach, we obtain
$$
L[s,t]=P[t]+O(s^2)+O(t^2),
$$ where $O(\cdot)$ is used in the sense of ideals in the ring of
formal powers series; the notation $+O(\cdot)$ means modulo the
addition of elements from the ideal generated by indicated
argument. Hence this implies $ \partial_t\partial_s L[0,0]= 0.  $
Since this is a contradiction we conclude that the operator in
\eqref{nop} cannot exist.

In what follows we use the notation
$$
D_{A_1\cdots A_k}:=D_{A_1}\cdots D_{A_k},\quad
\nd_{a_1\cdots a_k}:=\nd_{a_1}\cdots \nd_{a_k}.
$$
Let us write $ \nd^{(\ell)} C $ as a shorthand for the tensor 
$ \nd_{a_1\cdots a_\ell}C_{bcde}$, and set $\nd^{(0)}C=C$. 
Similarly we will write $ \nd^{(j)} \Delta^\ell f$ as a shorthand for
$\nd_{a_1\cdots a_j} \Delta^\ell f$.
Unless otherwise stated,
$\nd$, $C$, $\Rho$ are assumed to be defined with respect to $g_t$.
Finally we set 
$$
w=k-n/2,
$$
which is a positive integer. 
 
We will work with a one parameter family of metrics
$g_t$ such that
\begin{equation}\label{rho-cond}
  \Rho^{(\ell)}(p)=O(t^2), \qquad\ell\ge 2,
\end{equation}  
\begin{equation}\label{metric-cond1}
  \nd^{(\ell)}C(p)=O(t^2), \qquad 0\leq \ell\leq w+n-5,
\end{equation}  
\begin{equation}\label{metric-cond2}
  \nd^{(\ell)}\Delta^{n/2-2}\nd_{bc}
  C_{a}{}^b{}_{d}{}^c(p)=O(t^2),   
  \qquad  \ell\ge 0,
\end{equation} 
and, for $n\ge6$, 
\begin{equation}\label{metric-cond3}
 \Delta^{n/2-3}\nd_{b\,c\,(a_{w+2}\cdots a_3}C_{a_2}{}^b{}_{a_1)}{}^c(p)=O(t^2).
\end{equation}
We next take a scaling function $\Up$  such that
\begin{equation}\label{scale-cond}
 \nd_{(a_\ell\cdots a_1)}\Up(p,t)= 0,\quad \ell\ge2.
\end{equation} 
Finally, for the density $f$, we impose
\begin{equation} \label{sectionconds-nd}
\nd^{(\ell)}f  (p)=0,  \quad 0\leq \ell\leq w .
\end{equation}

From the conformal invariance of the Weyl curvature $C$ and
\eqref{transform} it follows that the condition \eqref{metric-cond1}
is conformally invariant --- in the sense that if a family of metrics
$g_t$ satisfies a set of conditions then so does $e^{2\Up}g_t$ for any
scaling function $\Up(x,t)$.  For $g_t$ satisfying 
\eqref{metric-cond1}, it is also
clear from \eqref{transform} that \eqref{metric-cond3}, which is a
condition on $\nd^{(w+n-4)}C$, is conformally invariant.  The
condition \eqref{metric-cond2} can be rewritten in terms of the
conformally invariant Fefferman-Graham obstruction tensor $B_{ab}$,
$$
  \nd^{(\ell)}B(p)=O(t^2),\quad   \ell\ge 0,
$$
because of \eqref{bach}.
Hence it is also conformally invariant.  
The condition \eqref{sectionconds-nd} is exactly equivalent to requiring
that $f$ is a density such that its $w$-jet vanish at $p$. Thus this
condition is independent of the choice of
$g_t$ and, in particular, is conformally invariant.  
Finally from \eqref{Rhotrans2} it is clear that
\eqref{rho-cond} is not a conformally invariant condition. Whereas the
point of conditions \eqref{metric-cond1},\eqref{metric-cond2} and
\eqref{metric-cond3} is to specialise the class of conformal structures
we allow, the role of \eqref{rho-cond} is rather as a scale normalisation
condition which restricts metrics allowed from within a given
conformal class $[g_t]$.  Nevertheless it is crucial to our arguments
that with \eqref{rho-cond} some conformal scaling freedom remains. In
particular if we assume \eqref{scale-cond} then we have \eqref{Rho-st}
below.

Under these assumptions, we will show the following two results.
\begin{lemma}\label{lemma-zero}
Assume that  \eqref{rho-cond}--\eqref{sectionconds-nd} hold.
Then 
$$
L[s,t]=P[t]+O(s^2)+O(t^2).
$$
\end{lemma}
\begin{lemma} \label{lemma-nzero}
Assume that \eqref{rho-cond}--\eqref{sectionconds-nd} hold
and further assume
\begin{equation} \label{sectionconds-harmonic}
\nd^{(\ell)}\Delta f (p)=O(t),  \qquad \ell\ge0.
\end{equation} 
Then
$$
  \partial_s L[0,t]=c\,\Up_b 
  \big(\Delta^{n/2-2}\nd_{c\,a_{w+1} \cdots\, a_{3}}
  C_{a_2}{}^b{}_{a_1}{}^c\big)
  \nd^{a_{w+1}\cdots\, a_{1}}f(p)
 +O(t^2),
$$
where $c$ is a non-zero constant.
\end{lemma}
With these
lemmas established the theorem is a consequence of the
existence of a density satisfying \eqref{sectionconds-harmonic} and the
following proposition which will be proved in the next section.

\begin{proposition}\label{metricprop}
There is a deformation $\{g_t\}_{t\in\bR}$ of the flat metric $g_0$
that satisfies \eqref{rho-cond}--\eqref{metric-cond3} yet with
\begin{equation}\label{dt-nonzero}
 F^b{}_{(a_{w+1}\cdots \,a_1)_0}(p)\ne0,
\end{equation}
where 
$$
 F^b{}_{a_{w+1}\cdots a_{1}}=
 \partial_t\big|_{t=0}\Delta^{n/2-2}
 \nd_{c\,a_{w+1} \cdots a_{3}}
   C_{a_2}{}^b{}_{a_{1}}{}^c.
$$
\end{proposition}

\noindent
{\bf Proof of Theorem
:}
Let $g_t$ be a family of metrics satisfying the
conditions \eqref{rho-cond}--\eqref{metric-cond3} and 
\eqref{dt-nonzero}. The existence of such a family is guaranteed by 
Proposition \ref{metricprop} above. 
Then by \eqref{dt-nonzero} we may find $\mu_b\in\ce_b|_p$ and
$
  \xi_{a_1a_2\cdots a_{w+1}}
  \in\ce_{(a_1\cdots a_{w+1})_0}[w]|_p
$
such that 
\begin{equation}\label{nonzero-linear}
 \mu_b \xi^{a_1a_2\cdots a_{w+1}}F^b{}_{a_1\cdots a_{w+1}}\ne0.
\end{equation}
Denoting by $x^i$ some fixed choice of normal coordinates for $g_0$ 
centered at $p$, we set
$$
f(x)=\xi_{i_1i_2\cdots i_{w+1}}x^{i_1}\cdots x^{i_{w+1}}.
$$ 
Then $f$ clearly satisfies the assumptions \eqref{sectionconds-nd}. In
the metric $g_0$ we also have $\nd^{(\ell)}\Delta f (p)=0$ for all $
\ell\ge0$. Thus, since the contorsion tensor distinguishing the metric
connections of $g_t$ and $g_0$ is $O(t)$, it follows immediately that $f$ 
satisfies \eqref{sectionconds-harmonic}.  Next we construct $\Up(x,t)$ 
by setting
$$
\Up(x,0)=\mu_ix^i
$$ 
and then obtain a function on $M\times \bR$ satisfying
\eqref{scale-cond} by solving the equation $\nd_{(a_1\cdots
a_\ell)}\Up(0,t)=0$, $\ell\ge2$, for each $t$.  From standard theory
this can be achieved within $ C^\infty(M\times \bR)$.
Then Lemma \ref{lemma-zero} implies
$L[s,t]=P[t]+O(s^2)+O(t^2)$ so that $\partial_s\partial_t L[0,0]=0$,
while Lemma \ref{lemma-nzero} shows
$$
 \partial_s\partial_t L[0,0]=c\, \mu_b \xi^{a_1a_2\cdots a_{w+1}}
 F^b{}_{(a_1\cdots a_{w+1})_0}\ne0,
$$
which is a contradiction.
\qed

\medskip

We now prove the lemmas used in the proof above.  Throughout the
proofs we will work at $ p\in M$.  In all final expressions the
tensors involved are evaluated at $p$ and we write simply $\nd_a
C_{bcde}$ to mean $\nd_a C_{bcde}(p)$ and so forth.

\medskip

\noindent{\bf Proof of Lemma \ref{lemma-zero}:}
We first prove 
$$
 S[t]:=L[0,t]-P[t]=O(t^2).
$$
Since  $L_g$ and $(-1)^kP_{2k}$ are natural operators which
agree for conformally flat metrics it follows that there is an 
expression for $S[t]$ as a sum of terms where 
each term is homogeneous of degree at least one in the jets of the 
curvature $R$ at $p$. Next via \eqref{R-PC} and \eqref{Bian} we may 
express the jets of $R$ in terms of $\Rho^{(\ell)}$ and
the jets of $C$ and obtain a new expression for $S[t]$ which is 
polynomial in these tensors. By \eqref{rho-cond} and since $C=O(t)$, 
the terms containing $\Rho^{(\ell)}$ and those which are nonlinear 
in $C$ are $O(t^2)$ and so can be neglected.
Thus using standard classical invariant theory and
elementary weight considerations we can express
the Riemannian invariant $S[t]$ mod $O(t^2)$ as a linear combination
of complete contractions of
$(\nd^{(\ell)}C)\nd^{(m)}f$ with $\ell+m=2k-2$.
In view of the conditions \eqref{metric-cond1} and 
\eqref{sectionconds-nd}, to obtain a non-trivial term
we must have $\ell\ge w+n-4$ and $m\ge w+1$. 
Thus a possible non-vanishing term should be a complete contraction of
one of the following two tensors:
$$
  \big(\nd^{(w+n-3)}C\big)\nd^{(w+1)} f\quad \text{ or } \quad
  \big(\nd^{(w+n-4)}C\big)\nd^{(w+2)} f . 
$$ 
Now consider the possible ways such a complete contraction could be
made.  First note that since the tensor field $C$ is completely
trace-free it is clear that in such a complete contraction we can assume, 
without loss of generality, that the indices of $C_{abcd}$ are paired 
with indices on $\nabla$'s. Next observe that the tensor field $C$
has the symmetry $ C_{abcd}=C_{[ab][cd]}$ while
$\nd^{a_{\ell}\cdots a_1}f=\nd^{(a_{\ell}\cdots a_1)}f+O(t)$.  
Thus $(\nd^{(\ell)}C_{abcd})\nd^{(m)}\nd^{abc}f=O(t^2)$ for
$\ell,$ $m\ge0$.  Also $\nd^{(\ell)}\nd^{abc}C_{abcd}=0$ and for both
similar results hold for any permutation of the indices on $C$. Thus
from the symmetries of the Weyl tensor $C$, the possible non-zero complete
contractions of the displayed terms must be contractions of the
tensors
$$
  \big(\nd^{(w-1)}\Delta^{n/2-2}\nd^{ab}C_{acbd}\big)
  \nd^{(w-1)}\nd^{cd}f
$$ 
and,  
$$
  \big(\nd^{(w)}\Delta^{n/2-3}\nd^{ab}C_{acbd}\big)
  \nd^{(w)}\nd^{cd}f, \quad \mbox{ if } n\geq 6 ,
$$
 or  
$$
\big(\nd^{(w-2)}\nd^{ab}C_{acbd}\big)
  \nd^{(w-2)}\Delta \nd^{cd}f, \quad \mbox{ if }  n=4 \mbox{ and } w\geq 2.
$$
But these are $O(t^2)$ by \eqref{metric-cond2} and 
\eqref{metric-cond3}. Thus $S[t]=O(t^2)$.  

To prove the general case, we first consider the tensors
$\Rho^{(\ell)}$ in the metric $g[s,t]$.
By the conformal transformation law \eqref{Rhotrans2} of $P^{(\ell)}$, we have
$$
\big[\Rho_{a_\ell\cdots a_1}\big]_{g[s,t]}
=\Rho_{a_\ell\cdots a_1}
-s\nd_{(a_\ell\cdots a_1)}\Up+O(s^2).
$$ 
Thus \eqref{rho-cond} and \eqref{scale-cond} imply 
\begin{equation}\label{Rho-st}
\big[\Rho_{a_\ell\cdots a_1}\big]_{g[s,t]}=O(s^2)+O(t^2).
\end{equation}
The other conditions imposed on, and properties of, $C$, $f$ and their
covariant derivatives, as used in the argument above at a metric $g_t$, are all
conformally invariant and so hold for $g[s,t]$. 
Thus replacing $g_t$ with $g[s,t]$ the argument
above that led to the conclusion $S[t]=O(t^2)$ can be repeated exactly
with the single exception that $\Rho^{(\ell)}$ can now be neglected
with error $O(s^2)+O(t^2) $ (rather than $O(t^2) $ as above).
Thus  with $ S[s,t]:=L[s,t]-P[t] $, we
obtain $S[s,t]=O(s^2)+O(t^2)$.  \qed

\medskip

\noindent{\bf Proof of Lemma \ref{lemma-nzero}:}
Since $D_{A_k\cdots A_1}$ is conformally invariant,
the conformal variation of $L_gf$ is determined entirely by the variation 
\eqref{XYZtrans} of $Y^{A}$. Thus we have 
$$
\partial_sL[0,t]=
 -\Upsilon_b\sum_{j=1}^{k}Z^{A_j b}Y^{A_{k}\cdots \widehat{A_{j}}
 \cdots A_1} D_{A_{k}\cdots A_{1}} f,
$$
where $\widehat{A_{j}}$ indicates an absent index.
Now (with $s$ still set to 0) we work with the metrics $g_t$.
From the definition of $D_{A_k\cdots A_1}$ we may re-express it as follows,
$$ 
\begin{aligned}
   D_{A_k\cdots A_1}f=&
  2D_{A_k\cdots[A_{j+1} A_{j}]\cdots A_1}f+
  2D_{A_k\cdots [A_{j+2} A_{j}]A_{j+1}\widehat{A_j}\cdots A_1}f+ \cdots\\ &+
  2D_{[A_k A_j]A_{k-1}\cdots \widehat{A_j}\cdots A_{1}}f
  +D_{A_jA_k\cdots \widehat{A_j}\cdots A_1}f .
\end{aligned}
$$
Using this we have at once that, 
$$
  \partial_sL[0,t]=-\sum_{j=1}^{k}j\Upsilon_bF^b_{(j)},
$$
where
$$
  F_{(j)}^b:=2 Z^{Bb}Y^{A_{k} \cdots\widehat{A_j} \cdots A_1}
               D_{A_k\cdots [A_{j+1}B]A_{j-1}\cdots A_1}f
$$
for $j\le k-1$, and 
$$
  F_{(k)}^b:= Z^{ Bb}Y^{A_{k-1} \cdots A_1}D_{BA_{k-1}\cdots A_1}f.
$$

We first show that $F_{(j)}^b=O(t^2)$ if $j\ne2$. Note that
$F_{(1)}^b=0$ because $D_{[A_2B]}$ vanishes on densities. Next we
recall that $f$ is a density of weight $k-n/2$ and each $D$ lowers
weight by 1. So from \eqref{Dform} we have
$$
Z^{Bb}D_{BA_{k-1}\cdots A_1}f=-Z^{Bb}X_B\Box D_{A_{k-1}\cdots A_1}f=0,
$$
which implies $F_{(k)}^b=0$.
To prove the other cases, we recall \eqref{Dcomm}:
\begin{equation}\label{commutator}
  D_{[CB]A_{j-1}\cdots A_{1}} f =(k-j)
  \sum_{\ell=1}^{j-1}E_{BC A_\ell}{}^Q D_{A_{j-1}
  \cdots  A_{\ell+1}QA_{\ell-1}\cdots A_{1}}f,
\end{equation}
where $E_{ABCD}$ is given by \eqref{def-E} and it is $O(t)$.
If we commute the indices for $D_{A_{j-1}\cdots A_{\ell+1}QA_{\ell-1}\cdots A_{1}}f$,
we get another $O(t)$ term.  Thus we see that each summand of the right-hand side of 
\eqref{commutator} is independent of $\ell$ up to permutations of $
A_1\cdots A_{j-1}$ and modulo $O(t^2)$.  
Hence  $F_{(j)}^b$ mod $O(t^2)$ for $j>2$ is a multiple of
\begin{equation}\label{F-zero}
Z^{Bb}Y^{A_{k-1}  \cdots A_1}
               D_{A_{k-1}\cdots A_{j+1}} (E_{A_{j}BA_{j-1}}{}^Q 
               D_{QA_{j-2}\cdots A_1}f ).
\end{equation}
This is $O(t^2)$, which we see as follows. From 
\begin{equation}\label{YDformula}
   [Y^A,\nd_b]=O(t)
\end{equation}
and the formula \eqref{Dform} for $D_B$, we conclude that there is a 
(weight dependent) operator $E_B$ such that 
$Y^AD_B=E_B Y^A+O(t)$ and so in \eqref{F-zero} we may pass  $Y^{A_1}$ to
the right where we finally observe that 
$$
Y^{A_1}D_{A_1}f=-\Delta f+O(t)=O(t),
$$
from \eqref{sectionconds-harmonic}.

We now focus on the computation of 
$$
 F_{(2)}^b=2 Z^{Bb}Y^{A_{k-1}\cdots A_1}
               D_{A_{k-1}\cdots A_3} ( E_{A_2BA_1}{}^Q 
               D_{Q}f ).
$$
Using \eqref{YDformula} and that $Y^AD_A=-\Delta+O(t)$,
we simplify $F_{(2)}^b$ to
\begin{equation}\label{typeII}
  2Z^{Bb}(-\Delta)^{k-3} ( Y^{A_2A_1}E_{A_2B A_{1}}{}^Q
  D_Qf ),
\end{equation}
modulo $O(t^2)$. Substituting \eqref{YYE},
 we may reduce the above formula to
a non-zero multiple of
\begin{equation}\label{F2-2}
Z^{Bb}\Delta^{k-3} \Big( Z_{B}{}^a Z^{Qq}\big(
(\nd^{cd}C_{cqda})
+2(\nd^cC_{cqda}) \nd^{d}\big)D_Qf \Big).
\end{equation}
Next using the identities of \eqref{connids} and figure \ref{XYZfigure},
we have
$$
Z^{Bb}\Delta^{k-3}Z_{B}{}^a=\bg^{ab}\Delta^{k-3}+O(t).
$$
Similarly, \eqref{connids} and \eqref{Dform} imply
$$
 (\nd^dZ^{Qq})D_Qf=-\bg^{dq}Y^QD_Qf=\bg^{dq}\Delta f
 +O(t)=O(t)
$$
and
$$
 Z^{Qq} D_Qf=(n+2w-2)\nd^{q}f.
$$
Thus \eqref{F2-2} is reduced, up to a non-zero multiple and modulo $O(t^2)$, to
\begin{equation}\label{reduction1}
\Delta^{k-3} \big( ( \nd^{cd}C_{cqd}{}^b) \nd^q f \big)
+2\Delta^{k-3}\big( ( \nd^cC_{cqd}{}^b)\nd^{dq} f\big).
\end{equation}
Finally we expand $\Delta^{k-3}$ by using the Leibniz rule.  Observe
that, in each term of the result, the total number of $\nd$ is
$2k-3=n+2w-3$, while in order to get a non-vanishing term, we need to
apply at least $(n+w-4)$ $\nd$'s to $C$, by \eqref{metric-cond1}, and
at least $(w+1)$ $\nd$'s to $f$, by \eqref{sectionconds-nd}. Such a
partition of $n+2w-3$ is unique and, using
\eqref{sectionconds-harmonic}, we see \eqref{reduction1} is reduced to
$$
2^w\begin{pmatrix}k-3\\ w\end{pmatrix}I^b
+2^{w}\begin{pmatrix}k-3\\ w-1\end{pmatrix} J^b,
$$ 
where
$$
\begin{aligned}
I^b&=
\big(\Delta^{n/2-3}\nd_{c\,d\,a_{w+1}\cdots a_2}C^c{}_{a_1}{}^{db}\big) 
  \nd^{a_{w+1}\cdots a_1} f,\\
J^b&=  
  \big(\Delta^{n/2-2}\nd_{c\,a_{w+1}\cdots a_3} C^c{}_{a_2a_1}{}^b\big) 
  \nd^{a_{w+1}\cdots a_1} f.
\end{aligned}
$$
If $n=4$, the first term does not appear 
and 
we immediately see that
$\partial_sL[0,t]$ mod $O(t^2)$ is
a non-zero multiple of ${J}^b\Up_b$.
If $n\ge6$, we have
\begin{equation}\label{relation-IJ}
(w+2)I^b=w J^b+ O(t^2),
\end{equation}
and, since $w>0$, are led to the same conclusion.
Therefore it remains only to prove this equation for $n\ge6$.

Corresponding to the curvature terms of $I^b$ and $J^b$, we set
$$
\begin{aligned}
{\widetilde I}&=\Delta^{n/2-3}\nd_{cd(a_{w+1}\cdots a_2}C^c{}_{a_1)}{}^{d}{}_{b},
\\
{\widetilde J}&=
 \Delta^{n/2-2}\nd_{c(a_{w+1}\cdots a_3} C^c{}_{a_2a_1)b},
\end{aligned}
$$ 
where we have suppressed the indices on ${\widetilde I} $ and
${\widetilde J}$ to simplify the notation.  Then noting that
$\nd_{cda_{w+1}\cdots a_2}C^c{}_{a_1}{}^{d}{}_b =\nd_{cd(a_{w+1}\cdots
a_2)}C^c{}_{(a_1}{}^{d}{}_{b)} +O(t^2)$, we may rewrite
\eqref{metric-cond3} as
\begin{equation}\label{IK-relation}
2{\widetilde I}+wK=O(t^2),
\end{equation}
where
$$
K=\Delta^{n/2-3}\nd_{cdb(a_{w+1}\cdots a_3}C^c{}_{a_2}{}^{d}{}_{a_1)}.
$$
On the other hand, note that from \eqref{Bian},
$$
\begin{aligned}
\nd_{[b} \nd^d C_{ca_2]d{a_1}} &=
2(n-3) \nd_{[b}\nd_c P_{a_2]a_1}\\
&= O(t^2),
\end{aligned}
$$
and hence $\nd^{cd}{}_{[b} C_{ca_2]d{a_1}} =O(t^2)$ which further implies
$$
\Delta^{n/2-3}\nd_{a_{w+1}\cdots a_3}\nd^{cd}{}_{[b} C_{ca_2]d{a_1}} = O(t^2).
$$
Symmetrising the left-hand-side of this last expression 
over $a_{w+1},\dots,a_1$ gives
$$
-\widetilde I+\widetilde J+K=O(t^2).
$$
Comparing this equation with \eqref{IK-relation},
we finally get
$$
(w+2)\widetilde I=w\widetilde J+O(t^2),
$$
which implies \eqref{relation-IJ}
because $ \nd^{a_{w+1}\cdots a_1} f= \nd^{(a_{w+1}\cdots a_1)} f+O(t)$.
\qed

\section{Construction of the metric} \label{metric}

It is clear that the issue of existence/nonexistence of invariant
operators is independent of signature (and could equally be treated in
the complex setting). To simplify the proof below we shall be
satisfied with constructing a Riemannian signature metric. With very
slight modification the same argument yields a proof of Proposition
\ref{metricprop} in any other desired signature.

\medskip

\noindent
{\bf Proof of Proposition \ref{metricprop}:}
 We first linearise the
problem.  For a symmetric two form $\psi=\psi_{ab}\in\ce_{(ab)}$ 
and each $t\in \bR$, we write $R_{abcd}[t], C_{abcd}[t]$ and 
$P_{ab}[t]$ respectively for the Riemannian curvature, the Weyl 
curvature and the Schouten tensor of $g_t=g_0+t\psi$.  Then set
$$
  \RP_{abcd}:=\frac{d}{dt}\Big|_{t=0}R_{abcd}[t], \quad 
  \CP_{abcd}:=\frac{d}{dt}\Big|_{t=0}C_{abcd}[t], \quad 
\RhoP_{ab}:= \frac{d}{dt}\Big|_{t=0}\Rho_{ab}[t]. 
$$
It follows from the definition of curvature that
\begin{equation}\label{Rform}
  \RP_{abcd}=
  \frac{1}{2}\left(\nd_{c[a}\psi_{b]d}-\nd_{d[a}\psi_{b]c}\right).
\end{equation}
Then $\CP_{abcd}$ is the trace-free part of
this, while $\RhoP_{ab}$ is a scaled trace adjustment of a single trace
of \eqref{Rform}.  Here $\nd$ is defined with respect to the flat metric
$g_0$. In terms of these tensors, Proposition \ref{metricprop} 
is reduced to the following lemma. 

\begin{lemma}\label{metric-lemma}
For each $w\in \bN$ there exists a symmetric two form
$\psi_{ab}\in\ce_{(ab)}$ on $\bR^n$ such that
\begin{equation}\label{metric-cond0lin}
\nd_{(a_\ell \cdots a_3}\RhoP_{a_2a_1)}(0)=0, \quad\ell\geq 2,
\end{equation} 
\begin{equation}\label{metric-cond1lin}
\nd^{(\ell)}\CP(0)=0, \quad 0\leq \ell\leq w+n-5,
\end{equation}  
\begin{equation}\label{metric-cond2lin}
 \nabla^{(\ell)}\Delta^{n/2-2}\nd_{bc}
 \CP_{a}{}^b{}_{d}{}^c(0)=0,
\quad \ell\ge 0,
\end{equation} 
and, for $n\ge6$
\begin{equation}\label{metric-cond3lin}
 \Delta^{n/2-3}\nd_{bc(a_{w+2} \cdots a_3}
  \CP_{a_2}{}^b{}_{a_1)}{}^c(0)=0
\end{equation}
yet with
\begin{equation}\label{metric-cond4lin}
\Delta^{n/2-2}\nd_{c(a_{w+1} \cdots a_{3}}
   \CP_{a_2}{}^b{}_{a_{1})_0}{}^c(0) \neq 0.
\end{equation}
\end{lemma}

Before we prove this we need some background on the representation
theory used.  The irreducible finite dimensional representations of
$\SL(n)$ can be classified by Young diagrams.  We use the notation
$(\ell_1,\ell_2,\dots , \ell_{n-1})$ to indicate the representation
corresponding to a Young diagram with rows (beginning from the top) of
length $\ell_1\geq \ell_2\geq \cdots \geq \ell_{n-1}\geq 0$. 
We identify $\SL(n)$ with its defining representation, and via this 
standard action of $\SL(n)$ on $\bR^n$ there are tensorial realisations
of these representations. For example a tensorial realisation of
$(\ell_1,\ell_2,\dots , \ell_{n-1})$ is given by
$Y(\ell_1,\ell_2,\dots , \ell_{n-1})$ which denotes the space of
(covariant) $\bR^n$ tensors, of rank $\sum \ell_i$, with the manifest
symmetries
$$
F_{a_1\cdots a_{\ell_1} b_1\cdots b_{\ell_2}\cdots 
 d_1\cdots d_{\ell_{n-1}}}=
F_{(a_1\cdots a_{\ell_1})( b_1\cdots b_{\ell_2})
\cdots (d_1\cdots d_{\ell_{n-1}})}
$$
and so called `hidden' symmetries which can be described as follows:
first a complete symmetrisation over any $\ell_1+1$ of the indices 
annihilates $F$; if we exclude the set $a_1\cdots a_{\ell_1}$,
then a complete symmetrisation over any $\ell_2+1$ of the remaining
indices annihilates $F$; if we exclude the sets $a_1\cdots a_{\ell_1}$
and $b_1\cdots b_{\ell_2}$ then a complete symmetrisation over any
$\ell_3+1$ of the remaining indices annihilates $F$ and so on. To
simplify the notation, we will omit terminal strings of zeros. Thus we
write $(\ell_1,\ell_2)$ as a shorthand for $(\ell_1,\ell_2,0,\dots
,0)$ and similarly $Y(\ell_1,\ell_2)$ for the described tensorial
realisation of this.

On the space of
tensors of rank $\sum\ell_i$, there are different  
projections onto a space $Y(\ell_1,\ell_2,\dots ,
\ell_{n-1})$ according to different orderings of the indices. (These
are easily described explicitly \cite{FH}). There are
identities between these projections but we do not need these
details. Any such projection will (also) be denoted by
$Y(\ell_1,\ell_2,\dots , \ell_{n-1})$ and is termed a Young
symmetriser.

The finite dimensional $\SO(n)$-representations (where as
usual $n$ is even) are also classified by strings of integers, in this
case just $n/2$ of these, $[\ell_1,\dots ,\ell_{n/2}]$, where
$\ell_1\geq \cdots \geq \ell_{n/2-1} \geq |\ell_{n/2}|$ and if 
$n/2$ is odd then $\ell_{n/2}\geq 0$.  We omit
terminal strings of zeros in this case too.  Via the defining
representation, where we view $\SO(n)$ as the subgroup of
$\SL(n)$ preserving the standard metric $\delta$, these also have tensorial
realisations: $Y_0[\ell_1,\dots ,\ell_{n/2-1}]$ (i.e.,
$Y_0[\ell_1,\dots ,\ell_{n/2-1},0]$) is the subspace of
$Y(\ell_1,\dots ,\ell_{n/2-1})$ consisting of completely trace-free
tensors. Continuing this notation if, $\ell_{n/2}> 0$ then the
subspace of $Y(\ell_1,\dots ,\ell_{n/2})$ consisting of completely
trace-free tensors, will be denoted $Y_0[\ell_1,\dots ,\ell_{n/2}]$.
This is an irreducible O$(n)$-module but upon restriction to
$\SO(n)$ is either irreducible or further decomposes 
depending on the parity of $n/2$: if $n/2$ is odd then 
$Y_0[\ell_1,\dots ,\ell_{n/2}]$ is irreducible, while if $n/2$ is even 
 $Y_0[\ell_1,\dots ,\ell_{n/2}]$ decomposes into a direct sum 
of irreducible
representations, each an eigenspace of an action of the volume form. 
The latter are realisations of representations usually denoted
$[\ell_1,\dots ,\ell_{n/2}]$ and $[\ell_1,\dots ,-\ell_{n/2}] $.
(In the corresponding complex theory of $\SO(n,\bC)$ one obtains 
such a decomposition regardless of the parity of $n/2$.)
\medskip

\noindent
{\bf Proof of Lemma \ref{metric-lemma}:}
 Consider the irreducible $\SL(n)$
representation $(n+w-2,2)$:
$$
\begin{picture}(60,18)(0,0)
\put(0,20){\line(1,0){60}}
\put(0,10){\line(1,0){60}}
\put(0,0){\line(1,0){20}}
\put(0,0){\line(0,1){20}}
\put(10,0){\line(0,1){20}}
\put(20,0){\line(0,1){20}}
\put(30,10){\line(0,1){10}}
\put(33,12){\mbox{$\cdots$}}
\put(50,10){\line(0,1){10}}
\put(60,10){\line(0,1){10}}
\end{picture}
$$
  Viewed as a representation of $\SO(n)$
 by restriction this decomposes into a direct sum of $\SO(n)$
irreducible representations. Using elementary representation theory
\cite{FH} it is easily verified that, since $n$ is even, 
the representation $[w+1,1]$: 
$$
\begin{picture}(60,18)(0,0)
\put(0,20){\line(1,0){50}}
\put(0,10){\line(1,0){50}}
\put(0,0){\line(1,0){10}}
\put(0,0){\line(0,1){20}}
\put(10,0){\line(0,1){20}}
\put(20,10){\line(0,1){10}}
\put(23,12){\mbox{$\cdots$}}
\put(40,10){\line(0,1){10}}
\put(50,10){\line(0,1){10}}
\put(52,9){\mbox{${}_0$}}
\end{picture}
$$
(and in dimension four $[w+1,1]\oplus [w+1,-1]$) occurs exactly once 
as a summand in this decomposition. 

At the level of tensor realisations this means that $Y_0[w+1,1]$ is a
summand in the orthogonal decomposition of $Y(n+w-2,2)$ and so for  
$K$  any non-zero tensor from the space $Y_0[w+1,1]$
there
is a Young symmetriser $Y(n+w-2,2)$ which does not annihilate
$$
K \otimes
\underbrace{\delta \otimes \cdots \otimes \delta}_{n/2-1}.
$$ 
Applying this Young symmetriser to the displayed tensor, let us
denote the image by $\Psi$. Then $\Psi\in Y(n+w-2,2)$ and so letting
$m:=n+w-2$ we have, in particular, that $\Psi=\Psi_{a_1\cdots a_m b_1 b_2}
=\Psi_{(a_1\cdots a_m)( b_1 b_2)}$.

We can view $\Psi$ as a (constant) covariant tensor on $\bR^n$ as
an affine space and in this setting we define
$\psi_{j_1j_2}:=\Psi_{i_1\cdots i_m j_1 j_2}x^{i_1} \cdots x^{i_{m}}$,
where $x^i$ are the standard coordinates.  Let $p$ be the origin in
$\bR^n$ (with $n\geq 4$ even as usual) and take $g_0:= \sum_1^n
dx^i\cdot dx^i $ so the component matrix of $g_0$ is $\delta$. Then we
claim that $\psi$ is a solution to
\eqref{metric-cond0lin}--\eqref{metric-cond4lin}.

First note that,
since $\psi$ is homogeneous of degree $m$, 
\eqref{metric-cond0lin}--\eqref{metric-cond2lin} are satisfied except
possibly for $\ell=m$ in \eqref{metric-cond0lin} and $\ell=w-2$ in
\eqref{metric-cond2lin}. 
In both cases the tensors on the left-hand-sides are obtained by
algebraic operations of symmetrisation and tracing from the tensor
$\Psi$. These operations are $\SO(n)$-equivariant and, since also the 
map from $K$ to $\Psi$ is $\SO(n)$-equivariant, we 
see that the maps from $K$ to these tensors are 
$\SO(n)$-equivariant.
 The same comment applies to the
left-hand-side of \eqref{metric-cond3lin} (for $n\geq 6$). 
 Consider first \eqref{metric-cond0lin} with $\ell=m$.  Note
that $\nd_{(b_{m} \cdots b_3}\RhoP_{b_2b_1)}(0)$ takes values in
$Y(m)$ and that, as an $\SO(n) $-module, $(m)$ decomposes into a
direct sum $[m]\oplus [m-2]\oplus \cdots $ (terminating in $[1]$ or
$[0]$ according to the parity of $m$).  Thus if $\nd_{(b_{m} \cdots
  b_3}\RhoP_{b_2b_1)}(0)$ were non-zero for any $K$ then, by the
composition of the linear equivariant operation mentioned with the
projection to $\SO(n)$-irreducible components, this would imply, when
$n\geq 6$, the existence of a non-trivial $\SO(n)$-module homomorphism
$[w+1,1]\to [q]$ for some $q\in \bN$,
$$
\begin{picture}(60,20)(0,4)
\put(0,20){\line(1,0){50}}
\put(0,10){\line(1,0){50}}
\put(0,0){\line(1,0){10}}
\put(0,0){\line(0,1){20}}
\put(10,0){\line(0,1){20}}
\put(20,10){\line(0,1){10}}
\put(23,12){\mbox{$\cdots$}}
\put(40,10){\line(0,1){10}}
\put(50,10){\line(0,1){10}}
\put(52,9){\mbox{${}_0$}}
\end{picture} 
 \to \ 
\begin{picture}(50,20)(0,0)
\put(0,10){\line(1,0){40}}
\put(0,0){\line(1,0){40}}
\put(0,0){\line(0,1){10}}
\put(10,0){\line(0,1){10}}
\put(13,2){\mbox{$\cdots$}}
\put(30,0){\line(0,1){10}}
\put(40,0){\line(0,1){10}}
\put(42,-1){\mbox{${}_0$}}
\end{picture},
$$
or similarly in dimension 4 would imply the existence of a
non-trivial $\SO(n)$-module homomorphism $([w+1,1]\oplus [w+1,-1])\to
[q]$ for some $q\in \bN$.  These equivariant mappings are
impossible since for any $w\in \bN$, $[w+1,1]$ and $[q]$ (and
$[w+1,-1]$ if $n=4$) are distinct irreducible $\SO(n)$-modules. An
almost identical argument shows that \eqref{metric-cond3lin} holds.  As
an $\SO(n)$-space, $(w+2)$ decomposes orthogonally to $[w+2]\oplus
[w]\oplus \cdots$ and so if $\Delta^{n/2-3}\nd_{bc(a_{w+2} \cdots
  {a_3}}\CP_{a_2}{}^b{}_{a_1)}{}^c(0)$ were non-zero for any $K$ then
we would once again arrive at a contradiction by deducing a
non-trivial linear $\SO(n)$-mapping $[w+1,1]\to [q]$ (or
$([w+1,1]\oplus [w+1,-1])\to [q]$ if $n=4$) for some $q\in \bN$.
Finally we consider \eqref{metric-cond2lin} with $\ell=w-2\geq 0$. Note
first that $ \nabla^{(w-2)}\Delta^{n/2-2}\nd_{bc}
\CP_{a}{}^b{}_{d}{}^c(0) $ has values in $Y(w-2)\otimes Y(2)$. For the
$\SL(n)$ tensor product we have the decomposition $(w-2)\otimes (2)=
(w)\oplus (w-1,1)\oplus(w-2,2)$. The $\SO(n)$-branch components of
$(w)$, $(w-1,1)$ and $(w-2,2)$ are all modules of the form $[a,b]$
where $a+b\leq w$ (where $b$ may be 0). In particular $[w+1,1]$ 
(and $[w+1,-1]$ for the case $n=4$) are not 
summands in the $\SO(n)$-decomposition of $(w-2)\otimes (2)$ and so
arguing as in the previous cases we immediately conclude that $
\nabla^{(w-2)}\Delta^{n/2-2}\nd_{bc} \CP_{a}{}^b{}_{d}{}^c(0) $ must
vanish.

It remains to establish \eqref{metric-cond4lin}. First we observe that
there is a non-trivial $(n/2-1)$-fold trace of $\Psi$ which takes
values in $Y_0[w+1,1]$. Up to scale this inverts the map which inserts
$Y_0[w+1,1]$ as an orthogonal summand in $Y(m,2)$.  Next observe that
if we skew over the pairs $a_1b_1$ and $a_2b_2$ of $\Psi_{a_1a_2\cdots
a_m b_1 b_2}$ and then on the result symmetrise over the indices
$a_1\cdots a_m$ and also over the indices $ b_1b_2$ then the result is
a non-zero multiple of $\Psi$. In fact this composition of mappings
can be taken as (up to scale) the definition of the Young symmetriser
projection onto the space $Y(m,2)$ containing $\Psi$.  Since, up to a
non-zero scale, $\Psi_{a_1\cdots a_m b_1 b_2}$ is $\nd_{a_m\cdots a_1}
\psi_{b_1b_2} $ it follows from \eqref{Rform} that $\Psi$ is a symmetry
adjustment of $\nd^{(m-2)}\RP(0)$. On the other hand from
\eqref{metric-cond0lin} and the linearisations of \eqref{R-PC} and 
\eqref{Bian} it follows that $\nd^{(m-2)}\RP(0)$ is a trace adjustment of
$\nd^{(m-2)}\CP(0)$ (cf.\ the discussion of normal scale in section
\ref{tractorsect}). Combining these observations it follows that there
is a symmetry and trace operation on $\nd^{(m-2)}\CP(0)$ with a
non-trivial outcome taking values in $Y_0[w+1,1]$.  Using the
symmetries of $\CP$ and the linearisation of the Bianchi identity
\eqref{Bian-Weyl} it is easily established that, up to scale, 
\eqref{metric-cond4lin} is the unique possibility.  \qed

%
%


\begin{thebibliography}{99}

\bibitem{BEGo} T.N. Bailey, M.G. Eastwood, and A.R. Gover,
{\em Thomas's structure bundle for conformal, projective and related
structures}, Rocky\ Mountain\ J.\ Math., {\bf 24} (1994) 1191--1217.

\bibitem{Brsharp} T. Branson, {\em Sharp inequalities, the functional
    determinant, and the complementary series}, Trans.\ Amer.\ Math.\
    Soc., {\bf 347} (1995) 3671--3742.

\bibitem{Cap-Gover2} A.\ \v Cap and A.R. Gover, {\em Tractor bundles for
  irreducible parabolic geometries}, Global analysis and harmonic
  analysis (Marseille-Luminy, 1999) 129-154, S\'{e}min.\ Congr.,
  {\bf 4}, Soc.\ Math.\ France, Paris 2000.
  Preprint ESI 865,
  available for viewing on the internet at 
  {\tt http://www.esi.ac.at}

\bibitem{Cap-Gover} A.\ \v Cap and A.R. Gover, {\em Tractor calculi for
parabolic geometries}, Trans.\ Amer.\ Math.\ Soc., {\bf 354} (2002)
1511-1548. 

\bibitem{CY} S.-Y.A. Chang and P.\ Yang, {\em Partial differential
equations related to the Gauss-Bonnet-Chern integrand on $4$-manifolds},
in ``Conformal, Riemannian and Lagrangian Geometry: The 2000 Barrett Lectures,''
Amer.\ Math.\ Soc., 2002. 


\bibitem{ESlo} M.G. Eastwood and J. Slov\'ak,
\textit{Semiholonomic Verma modules}, J.\ Algebra, {\bf 197} (1997) 424--448.


\bibitem{FG2} C. Fefferman and C.R. Graham,
{\em Conformal invariants}, in {\em Elie Cartan et les math\-\'ematiques
d'aujourd'hui}, Ast\'erisque, hors s\'{e}rie (Soci\'{e}t\'{e}
Math\'{e}matique de France, Paris, 1985) 95--116.

\bibitem{FH} W. Fulton and J. Harris, Representation Theory:
A First course,  Springer, 1991.

\bibitem{Gosrni} A.R.\ Gover, {\em Aspects of parabolic invariant theory}, 
 Supp.\ Rend.\ Circ.\ Matem.\ Palermo, Ser.\ II, Suppl.\ {\bf 59}
(1999) 25--47.

\bibitem{Goadv} A.R.\ Gover, {\em Invariant theory and calculus for
    conformal geometries}, Adv.\ Math., {\bf 163} (2001) 206--257.

\bibitem{GoGr} A.R.\ Gover and C.R.\ Graham, {\em CR invariant powers of
    the sub-Laplacian},
preprint {\tt arXiv:math.DG/0301092}

\bibitem{GP} A.R.\ Gover and L.J.\ Peterson, 
{\em Conformally invariant powers of the Laplacian, Q-curvature, 
and tractor calculus},
Comm.\ Math.\ Phys., {\bf 235} (2003) 339--378.

\bibitem{Grno} C.R.\ Graham, {\em Conformally invariant powers of the
    Laplacian, II: Nonexistence}, J.\ London Math.\ Soc.\ (2) {\bf 46} (1992)
566--576.

\bibitem{GJMS} C.R.\ Graham, R. Jenne, L. Mason, and G. Sparling,
{\em Conformally invariant powers of the Laplacian, I: existence}, J.\
London Math.\ Soc.\ (2), {\bf 46} (1992) 557--565.

\bibitem{GZ} C.R.\ Graham and M.\ Zworski,
{\em Scattering matrix in conformal geometry},
Invent. math., {\bf 152} (2003) 89--118. 


\bibitem{JV} H.P. Jakobsen and M.\ Vergne, {\em Wave and Dirac
operators, and representations of the conformal group}, J.\ Func.\
Anal., {\bf 24} (1977) 52--106.

\bibitem{ot} R. Penrose and W. Rindler,  Spinors and Space-time, vol 1,
Cambridge Univ.\ Press, 1984.

\bibitem{T} T.Y.\ Thomas, {\em On conformal geometry},
  Proc.\ Natl.\ Acad.\ Sci.\ USA,
  {\bf 12} (1926) 352--359.

\bibitem{W} V. W\"unsch, 
\textit{Some new conformal covariants}, 
Z.\ Anal.\ Anw. {\bf 19} (2000) 339--357.
Erratum, Z.\ Anal.\ Anw., {\bf 21} (2002) 529--530.
   

\end{thebibliography}
\end{document}